\documentclass [a4,10pt,oneside]{article}
\usepackage{amsmath, amssymb, amsthm, amsfonts, amsxtra, latexsym, amscd,
pb-diagram,  graphics,rotating,xy}

\theoremstyle{plain}

\newcommand{\A}{\mathcal{A}}
\newcommand{\tx}{\otimes }
\newcommand{\ts}{\oplus}
\newcommand{\Lh}{\frak{L}}
\newcommand{\Rh}{\frak{R}}

\newcommand{\Fm}{\widetilde{F}}

\newcommand{\Fa}{\Breve{F}}

\newcommand{\Gm}{\widetilde{G}}
\newcommand{\Ga}{\Breve{G}}

\begin{document}

\title{ \Large{\bf Cohomological classification of Ann-functors}}
\author {Nguyen Tien Quang and Dang Dinh Hanh}

\maketitle


\begin{abstract}

Regular Ann-functor classification problem has been solved with Shukla cohomology. In this paper, we would like to present a solution to the above problem in the general case and in the case of strong Ann-functors with, respectively, Mac Lane cohomology and Hochschild cohomology.

\end{abstract}


\vspace{0.3cm}
 

\vspace{0.5cm}

\section{ Introduction}

The definition of  {\it Ann-categories} was presented by N.T.Quang [7] in 1987, which is regarded as a categorization of ring structure. Each Ann-category $\mathcal A$ is Ann-equivalent to its reduced Ann-category. This Ann-category is of the  type  $(R, M, h)$, where $R$ is a ring of congruence classes of objects of   $\mathcal A$, $M=Aut(0)$ is the  $R$-bimodule and $h$ is a 3-cocycle in $Z^3_{MaL}(R, M)$ (due to Mac Lane [6]). Then,  there exists a bijection between the set of congruence classes of Ann-categories of the type $(R,M)$ and the cohomology group $H^3_{MaL}(R, M)$ (see [11]). For the {\it regular Ann-categories} (whose the commutativity constraint satisfies $c_{X, X}=id$), then in the above bijection, the group  $H^3_{MaL}(R, M)$ is replaced with the Shukla cohomology group $H^3_{Sh}(R, M)$ [14].

In 2006 [4], M.Jibladze and T. Pirashvili  presented the definition of {\it categorical rings} as  a slightly modified version of the definition of Ann-categories and classified them by Mac Lane ring cohomology. However, in [10] authors have showed that, it has not been proved whether the $R$-bimodule structure on $M$ can be deduced from the axiomatics of categorical rings.

The  Ann-functor classification problem has been solved for regular  Ann-categories with Shukla cohomology  [1, 12]. In this paper, we  present a solution  for this problem in the general case via  low-dimensioned cohomology groups of  Mac Lane ring cohomology. In proper, Hochschild algebra cohomology used to classification of strong Ann-functor.

In this paper, for convenience, sometimes we denote by $XY$ the tensor product of the two objects $X$ and $Y$, instead of $X\otimes Y$.

\section{Preliminaries}
\subsection {The basic concepts}

The definition of Ann-categories was presented in [7, 10, 11]. We always suppose that  $\A$ is an  Ann-category with a system of constraints:
\[(a^+, c, (O, g, d), a, (I, l, r), \mathcal L, \mathcal R).\]

{\bf Definition 1.} {\it Let $\A$ and $\A'$ be Ann-categories. An Ann-functor from  
$\A$ to $\A'$ is a triple $(F, \Fa, \Fm)$, where $(F, \Fa)$ is a symmetric monoidal functor respect to the  operation $\oplus$, $(F, \Fm)$ is an $A$-functor (i.e. an associativity functor) respect to the operation  $\otimes$, satisfies  the two following commutative diagrams:

\[\scriptsize\begin{diagram}
\node{F(X(Y\oplus Z))}\arrow{s,l}{F(\Lh)}
\node{FX.F(Y\oplus Z)}\arrow{w,t}{\Fm}
\node{FX(FY\oplus FZ)}\arrow{w,t}{id\otimes \Fa}\arrow{s,r}{\Lh'}\\
\node{F(XY\oplus XZ)}
\node{F(XY)\oplus F(XZ)}\arrow{w,t}{\Fa}
\node{FX.FY\oplus FX.FZ}\arrow{w,t}{\Fm \oplus \Fm}
\end{diagram}\]
\[\scriptsize\begin{diagram}
\node{F((X\oplus Y)Z)}\arrow{s,l}{F(\Rh)}
\node{F(X\oplus Y).FZ}\arrow{w,t}{\Fm} 
\node{(FX\oplus FY).FZ}\arrow{w,t}{\Fa \otimes id}\arrow{s,r}{\Rh'}\\
\node{F(XZ\oplus YZ)}
\node{F(XZ)\oplus F(YZ)}\arrow{w,t}{\Fa}
\node{FX.FZ\oplus FY.FZ}\arrow{w,t}{\Fm \oplus \Fm}
\end{diagram}\]

The commutation of the  above diagrams are called the compatibility of the functor $F$ with the distributivity constraints of the two Ann-categories $\mathcal A,\mathcal A'$.

\indent We call $\varphi : F\rightarrow G$ an $Ann$-\emph{morphism} between two  $Ann$-functors $(F, \Fa, \Fm)$ and $(G, \Ga, \Gm)$ if it is an  $\oplus$-morphism as well as an $\otimes$-morphism.

An Ann-functor $(F,\Fa, \Fm):\A \rightarrow  \A'$ is called an \emph{Ann-equivalence} if\linebreak there  exists an Ann-functor  $(F',\Fa',\Fm'):\A'\rightarrow  \A$ and natural isomorphisms
\linebreak $\alpha: F\circ F' \cong id_{\A'},\quad \alpha' : F'\circ F \cong id_{\A}$. }

By  Theorem  8 [9], an Ann-functor $(F,\Fa, \Fm):\A \rightarrow  \A'$ { is an  Ann-equivalence iff   $F$ is a categorical equivalence.}

Note that, similar to a ring homomorphism, an Ann-functor $F$ is not required $F(1)\cong 1^{'}.$ Moreover, note that: in the Definition 1, it is only required that $(F, \Fa)$ is an $AC$-functor (i.e. an $\oplus$-functor which is compatible with the associativity and commutativity constraints). Indeed, since $(\A, \ts), (\A', \ts)$ are $Gr-$categories, each $A$-functor is compatible with the unitivity constraints.

\subsection{The third Mac Lane ring cohomology group $H^3_{MaL}(R,M)$}

Let $R$ be a ring and  $M$ be an $R$-bimodule. From the definition of  Mac Lane ring cohomology  [6], we may obtain the description of the elements of cohomology group $H^3_{MaL}(R,M)$.

The group  $Z^3_{MaL}(R,M)$ of  3-cochains of the ring $R$, with coefficients in\linebreak $R$-bimodules $M$, consisting of quadruples $(\sigma,\alpha,\lambda,\rho)$, functions:

\[\alpha,\lambda,\rho:R^3\rightarrow M\]
and  $\sigma:R^4\rightarrow M$ 
satisfy the following relations:\\

\noindent M1. $x\alpha(y,z,t)-\alpha(xy,z,t)+\alpha(x,yz,t)-\alpha(x,y,zt)+\alpha(x,y,z)t=0$

\vspace{0.1cm}
\noindent M2. $\alpha(x,z,t)+\alpha(y,z,t)-\alpha(x+y,z,t)+\rho(xz,yz,t)-\rho(x,y,zt)+\rho(x,y,z)t=0$

\vspace{0.1cm}
\noindent M3. $-\alpha (x, y,  t) - \alpha (x, z,t)+\alpha (x, y+z, t)
+  x\rho (y, z, t)-\rho (xy, xz, t)\\
\indent  - \lambda (x, yt, zt) +\lambda (x, y, z) t =0$

\vspace{0.1cm}
\noindent M4. $ \alpha(x,y,z)+\alpha(x,y,t)-\alpha(x,y,z+t)+x\lambda (y,z,t)-\lambda(xy,z,t)+\lambda(x,yz,yt)=0 $

\vspace{0.1cm}
\noindent M5. $\lambda(x,z,t)+\lambda(y,z,t)-\lambda(x+y,z,t)+\rho(x,y,z)+\rho(x,y,t)-\rho(x,y,z+t)\\
\indent +\sigma(xz,xt,yz,yt)=0$

\vspace{0.1cm}
\noindent M6. $\lambda(x,a,b)+\lambda(x,c,d)-\lambda(x,a+c,b+d)-\lambda(x,a,c)-\lambda(x,b,d)\\
\indent +\lambda(x,a+b,c+d)-x\sigma(a,b,c,d)+\sigma(xa,xb,xc,xd)=0$

\vspace{0.1cm}
\noindent M7. $\rho(a,b,x)+\rho(c,d,x)-\rho(a+c,b+d,x)-\rho(a,c,x)-\rho(b,d,x)\\
\indent +\rho(a+b,c+d,x)-\sigma(ax,bx,cx,dx)+\sigma(a,b,c,d)x=0$

\vspace{0.1cm}
\noindent M8. $\sigma (a,b,c,d)+\sigma (x,y,z,t)-\sigma (a+x,b+y,c+z,d+t)+\sigma (a,b,x,y)+\sigma (c,d,z,t)$\\
$-\sigma (a+c,b+d,x+z,y+t)
+\sigma (a,c,x,z)+\sigma (b,d,y,t)-\sigma (a+b,c+d,x+y,z+t)=0$

\vspace{0.1cm}
\noindent M9.  $\alpha (0,y,z)=\alpha(x,0,z)=\alpha(x,y,0)=0$

\vspace{0.1cm}
\noindent M10. $\sigma(0,0,z,t)=\sigma(x,y,0,0)=\sigma(0,y,0,t)=\sigma(x,0,z,0)=\sigma(x,0,0,t)=0$.\\

The subgroup $B^3_{MaL}(R,M)\subset Z^3_{MaL}(R,M)$ of 3-coboundaries consists of the quadruples $(\sigma,\alpha,\lambda,\rho)$ such that there exist the maps  $\mu,\nu:R^2\rightarrow M$ satisfying:\\

\noindent M11. \  $\sigma (x,y,z,t)=-\mu(x,y)-\mu(z,t)+\mu(x+z,y+t)+\mu(x,z)+\mu(y,t)$\\
\indent\qquad\qquad\qquad\qquad $-\mu(x+y,z+t)$

\vspace{0.1cm}
\noindent M12. \ \  \ $\alpha(x,y,z)=x\nu(y,z)-\nu(xy,z)+\nu (x,yz)-\nu(x,y)z$

\vspace{0.1cm}
\noindent M13. \ \  \ $\lambda(x,y,z)=\nu(x,y)+\nu(x,z)-\nu(x,y+z)+x\mu(y,z)-\mu(xy,xz)$

\vspace{0.1cm}
\noindent M14.\ \ \ \  $\rho(x,y,z)=-\nu(x,z)-\nu(y,z)+\nu(x+y,z)+\mu(xz,yz)-\mu(x,y)z.$\\

Finally,
$H^3_{MaL}(R,M)=Z^3_{MaL}(R,M)/B^3_{MaL}(R,M).$\\

Each Ann-category $\mathcal I$ of the type $(R, M)$ having the {\it structure   f} \  is a family $f=(\xi, \eta, \alpha, \lambda, \rho)$, where $\xi, \alpha, \lambda, \rho: R^3\rightarrow M$ and $\eta: R^2\rightarrow M$ are functions satisfying  17 the equations (see Proposition 5.8 [11]). Now, we define a function $\sigma: R^4\rightarrow  M$,  given by:
\[\sigma(x,y,z,t)=\xi(x+y,z,t)-\xi(x,y,z)+\eta(y,z)+\xi(x,z,y)-\xi(x+z,y,t).\]
This function is respect to the associativity-commutativity constraint $v$ in the Ann-category $\A$, where
$$v = v_{X,Y,Z,T}:(X\ts Y)\ts (Z\ts T) \longrightarrow (X\ts Z)\ts (Y\ts T)$$ 
is given by commutative diagram fowllowing:
\[\scriptsize\begin{diagram}
\node{(X\ts Y)\ts(Z\ts T)}\arrow{s,l}{v} \arrow{e,t}{a_{+}}
\node{((X\ts Y)\ts Z)\ts T}
\node{(X\ts (Y\ts Z))\ts T}\arrow{s,r}{(X\ts c)\ts T}\arrow{w,t}{a_+\ts T}\\
\node{(X\ts Z)\ts(Y\ts T)}\arrow{e,t}{a_{+}}
\node{((X\ts Z)\ts Y)\ts T}
\node{(X\ts(Z\ts Y))\ts T}\arrow{w,t}{a_+\ts T}
\end{diagram}\]
The quadruple $h=(\sigma,\alpha,\lambda,\rho)$ is a 3-cocycle of the ring $R$ with coefficients  in $R$-bimodule $M$ due to Mac Lane (Theorem 7.2 [11]) and therefore each reduced Ann-category is of the form $(R,M,h)$.

\section{ An equivalence criterion of an  Ann-functor}

First, we show a characterized property of Ann-functors, which is related to the associativity-commutativity constraint $v$.

{\bf Definition 2.} {\it Let $\A,\   \A'$ be symmetric  monoidal $\oplus-$categories. Then, the  $\oplus-$functor $(F,\Fa):\A\rightarrow \A'$ is called compatible with the constraints  $v, v^{'}$  if the following diagram commutes for all $X,Y,Z,T\in \A$}
{\scriptsize 
\[
\divide \dgARROWLENGTH by 2
\begin{diagram}
\node{F((X\ts Y)\ts(Z\ts T))}\arrow{s,l}{F(v)}
\node{F(X\ts Y)\ts F(Z\ts T)}\arrow{w,t}{\Fa}
\node{(FX\ts FY)\ts (FZ\ts FT)}\arrow{w,t}{\Fa+\Fa}\arrow{s,r}{v'}\\
\node{F((X\ts Z)\ts (Y\ts T))}
\node{F(X\ts Z)\ts F(Y\ts T)}\arrow{w,t}{\Fa}
\node{(FX\ts FZ)\ts (FY\ts FT)}\arrow{w,t}{\Fa+\Fa}\tag{1}
\end{diagram}
\]
}

\noindent Then 
{\bf Lemma 3.1.} {\it 
Let $\oplus-$functor  $(F,\Fa):\A\rightarrow \A'$ be compatible with the unitivity constraints. Then $(F,\Fa)$ is an $AC-$functor iff it is compatible with the constraints $v, v^{'}$.}
\begin{proof}

The nescessary condition was presented by D. B. A. Epstein (Lemma 1.5 [2]).

Now,  assume that the diagram (1) commutes. To prove that the pair  $(F,\Fa)$ is compatible with the commutativity constraints, we consider the following Diagram 1. 
 
In the Diagram 1,  the region (I) commutes thanks to the naturality of the morphism $v$,\ the regions (II) and (IV) commute since $(F,\Fa)$  is compatible with the unitivity constraints, the regions (III) and (VII) commute by the coherence theorem in a symmetric monoidal category, the regions  (VI) and (VIII) commute thanks to the naturality of  $\Fa$, the outside region commutes by the diagram (1). Hence, the region (V) commutes. So $(F,\Fa)$ is compatible with the commutativity constraints.

\[\scriptsize\setlength{\unitlength}{0.5cm}
\begin{picture}(26,17)

\put(0,1){$(FO\ts FX)\ts (FY\ts FO)$}
\put(16,1){$(FO\ts FY)\ts (FX\ts FO)$}

\put(5,4){$(O\ts FX)\ts (FY\ts O)$}
\put(13,4){$(O\ts FY)\ts (FX\ts O)$}

\put(8,7){$FX\ts FY$}
\put(13,7){$FY\ts FX$}

\put(8,13){$F(X\ts Y)$}
\put(13,13){$F(Y\ts X)$}

\put(0,10){$F(O\ts X)\ts F(Y\ts O)$}
\put(17,10){$F(O\ts Y)\ts F(X\ts O)$}

\put(0,16){$F((O\ts X)\ts (Y\ts O))$}
\put(17,16){$F((O\ts Y)\ts (X\ts O))$}

\put(7, 1.2){\vector(1,0){8.5}}\put(11.5,1.4){$v'$}
\put(3.5, 1.7){\vector(1,1){1.9}}\put(5.2,2.5){$(\widehat F\ts id)\ts (id\ts \widehat F)$}
\put(20.5, 1.7){\vector(-1,1){1.9}}\put(14,2.5){$(\widehat F\ts id)\ts (id\ts \widehat F)$}

\put(10.8, 4.2){\vector(1,0){2.05}}\put(11.5,4.4){$v'$}

\put(9.3, 4.6){\vector(0,1){2.1}}\put(7.1,5.5){$g'\ts d'$}

\put(5.5, 9.7){\vector(1,-1){2.3}}\put(3,8.5){$F(g)\ts F(d)$}

\put(6, 15.8){\vector(1,-1){2.1}}\put(4.5,14.5){$F(g\ts d)$}

\put(16.8, 9.7){\vector(-1,-1){2.1}}\put(16.4,8.5){$F(g)\ts F(d)$}
\put(16.5, 15.8){\vector(-1,-1){2.1}}\put(16,14.5){$F(g\ts d)$}

\put(9.3, 7.7){\vector(0,1){5.1}}\put(8.6,10){$\Fa$}

\put(14.2, 4.7){\vector(0,1){2}}\put(14.5,5.5){$g'\ts d'$}

\put(14.2, 7.7){\vector(0,1){5.1}}\put(14.5,10){$\Fa$}

\put(10.8, 7.2){\vector(1,0){2.1}}\put(11.5,7.5){$c'$}

\put(10.7, 13.2){\vector(1,0){2.1}}\put(11.1,12.5){$F(c)$}

\put(1.8, 1.7){\vector(0,1){8}}\put(0,5.5){$\Fa\ts \Fa$}
\put(1.8, 10.7){\vector(0,1){5}}\put(0.5,13){$\Fa$}

\put(21, 1.7){\vector(0,1){8}}\put(21.4,5.5){$\Fa\ts \Fa$}
\put(21, 10.7){\vector(0,1){5}}\put(21.4,13){$\Fa$}

\put(6, 16.1){\vector(1,0){10.5}}\put(11,16.3){$F(v)$}

\put(11.3,2.5){(I)}
\put(11.2,5.5){(III)}
\put(11.3,10){(V)}
\put(11.1,14.5){(VII)}

\put(4,5.5){(II)}
\put(4,13){(VI)}

\put(19,5.5){(IV)}
\put(19,13){(VIII)}
\put(10, 0){Diagram 1}
\end{picture}\]

Next, we consider the following Diagram 2.

In the Diagram 2, the region  (I) commutes thanks to the  naturality of the morphism $v$;
the first component of the region (II) commutes since $(F, \Fa)$ is compatible with unitivity constraints, the second one commutes thanks to the composition of morphisms, so  the region (II) commutes; the regions (III) and (X) commute thanks to the coherence in a symmetric monoidal category; the first component of the region (IV) commutes thanks to the composition of morphisms, the second one commutes since $(F,\Fa)$  is compatible with unitivity constraints, so the region (IV) commutes; the region (V) and (VII)  commute thanks to the composition of morphisms; the regions (VIII) and (IX)  commute thanks to the naturality of $\Fa$; the outside region commutes thanks to the diagram (1). Therefore, the region (V) commutes, i.e., the pair $(F,\Fa)$ is compatible with associtivity constraints.

\[\scriptsize
\setlength{\unitlength}{0.5cm}
\begin{picture}(26,19)

\put(0,1){$(FX\ts FO)\ts (FY\ts FZ)$}
\put(16,1){$(FX\ts FY)\ts (FO\ts FZ)$}

\put(3.5,4){$(FX\ts O)\ts (FY\ts FZ)$}
\put(13,4){$(FX\ts FY)\ts (O\ts FZ)$}

\put(5.5,7){$FX\ts (FY\ts FZ)$}
\put(13,7){$(FX\ts FY)\ts FZ$}

\put(0,10){$F(X\ts O)\ts F(Y\ts Z)$}
\put(17,10){$F(X\ts Y)\ts F(O\ts Z)$}

\put(5.5,13){$FX\ts F(Y\ts Z)$}
\put(13,13){$F(X\ts Y)\ts FZ$}

\put(5.75,16){$F(X\ts (Y\ts Z))$}
\put(12.7,16){$F((X\ts Y)\ts Z)$}

\put(0,19){$F((X\ts O)\ts (Y\ts Z))$}
\put(17,19){$F((X\ts Y)\ts (O\ts Z))$}

\put(6.7, 1.1){\vector(1,0){9}}\put(11.2,1.2){$v'$}
\put(3.5, 1.5){\vector(1,1){2.2}}\put(5.3,2.5){$(id\ts \widehat F)\ts id$}

\put(19, 1.5){\vector(-1,1){2.2}}\put(14,2.5){$id\ts (\widehat F\ts id)$}

\put(9.7, 4.2){\vector(1,0){3.1}}\put(11.2,4.4){$v'$}

\put(8.5, 4.65){\vector(0,1){2.1}}\put(6.6,5.5){$g'\ts id$}

\put(14.5, 4.65){\vector(0,1){2.1}}\put(14.9,5.5){$id\ts d'$}

\put(10.2, 7.2){\vector(1,0){2.6}}\put(11.2,7.5){$a_+$}

\put(7.7, 7.5){\vector(-1,1){2.1}}\put(3,8.5){$F(g^{-1})\ts \Fa$}
\put(14.9, 7.5){\vector(1,1){2.2}}\put(16.7,8.5){$\Fa\ts F(d^{-1})$}

\put(8.5, 7.5){\vector(0,1){5.2}}\put(8.8,8.5){$id\ts \Fa$}
\put(14.5, 7.5){\vector(0,1){5.2}}\put(12.5,11.5){$\Fa\ts id$}

\put(5.7, 10.5){\vector(1,1){2.1}}\put(3.9,11.5){$F(g)\ts id$}
\put(17.2, 10.6){\vector(-1,1){2}}\put(16.7,11.5){$id\ts F(d)$}

\put(8.5, 13.65){\vector(0,1){2.1}}\put(8.8,14.5){$\Fa$}
\put(14.5, 13.65){\vector(0,1){2.1}}\put(13.8,14.5){$\Fa$}

\put(10, 16.2){\vector(1,0){2.5}}\put(10.4,16.5){$F(a_+)$}

\put(16.6, 18.8){\vector(-1,-1){2.2}}\put(16.2,17.5){$F(id\ts d)$}
\put(6, 18.8){\vector(1,-1){2.2}}\put(4.2,17.5){$F(g\ts id)$}

\put(6, 19.1){\vector(1,0){10.6}}\put(10.6,19.3){$F(v)$}

\put(1.8, 1.7){\vector(0,1){8}}\put(0,5.5){$\Fa\ts \Fa$}
\put(1.8, 10.7){\vector(0,1){8}}\put(1,16){$\Fa$}

\put(21, 1.7){\vector(0,1){8}}\put(21.3,5.5){$\Fa\ts \Fa$}
\put(21, 10.7){\vector(0,1){8}}\put(21.3,16){$\Fa$}

\put(11.1,2.5){(I)}
\put(10.9,5.5){(III)}
\put(7,10){(V)}
\put(10.9,10){(VI)}
\put(15,10){(VII)}
\put(4,14.5){(VIII)}
\put(19,14.5){(IX)}

\put(10.9,17.7){(X)}

\put(4,5.5){(II)}

\put(19,5.5){(IV)}
\put(10, 0){Diagram 2}
\end{picture}\]

\end{proof}

{\bf Proposition 3.2.} {\it

In the definition of Ann-functors, the condition that $(F, \Fa)$ is an symmetric monoidal $\oplus$-functor is equivalent to the two following conditions:
\begin{enumerate}
\item $(F, \Fa)$ is compatible with the unitivity constraints respect to the operation $\oplus$,

\item $(F, \Fa)$ is compatible with the constraints  $v, v^{'}$.
\end{enumerate}
}
\begin{proof}
 Directly deduced from Lemma 3.1.
\end{proof}

\section{Ann-functors and the low-dimensioned cohomology  groups of rings due to  Mac Lane}
\subsection{  Ann-functors of the type $(p, q)$}

Now, we will show that each Ann-functor $(F, \Fa,\Fm):\A\rightarrow \A^{'}$ induces a Ann-functor
$\overline{F}$ on their reduced Ann-categories, and this correspondence is 1-1. First, we have

{\bf Theorem 4.1.}
[{Theorem 4.6 [11]}]
{\it Let $\A$ and $\A'$ be Ann-categories. Then, each Ann-functor $(F,\breve F, \Fm):\A\rightarrow \A'$ induces the pair of ring homomorphisms:
\[\begin{array}{ccccccccccc}
F_{0}:&\Pi_0(\A) &\rightarrow& \Pi_0(\A')&; \qquad \qquad &\mathop F_{1}:&\Pi_1(\A) &\rightarrow& \Pi_1(\A')\\
                      & cls X        &\mapsto      &  cls FX    &                            &                                          & u              & \mapsto     & \gamma^{-1}_{F0}(Fu)
\end{array}\]
satisfying 
$$
 F_{1}(su) =F_{0}(s) F_{1}(u); \qquad F_{1}(us) = F_{1}(u)F_{0}(s)
$$
where $\Pi_1(\A)$ is regarded as a ring with the null multiplication. Furthermore, $F$ is a equivalence iff $F_{0}, F_1$ are isomorphisms.
}

The pair  $(F_0,F_1)$ is called {\it the pair of induced homomorphisms} of the $Ann-$functor $(F,\Fa,\Fm)$. If $\mathcal S, \mathcal S'$ are, respectively,  the reduced  $Ann-$categories  of $\A, \A'$ then the functor $\overline{F}: \mathcal S\rightarrow \mathcal S'$ given by
$$\overline F(s)=F_0(s),\  \overline F (s,u)=(F_0 s,F_1 u)$$
is called  the {\it reduced functor} of  $(F,\Fa, \Fm)$ on reduced $Ann-$categories.

{\bf Proposition 4.2.} {\it
Let $\overline F$ be the induced functor of the   $Ann-$functor  $(F,\Fa, \Fm): \A\rightarrow \A'$. Then the diagram
\[\begin{diagram}
\node{\A}\arrow{e,t}{F}
\node{\A'}\arrow{s,r}{G'}\\
\node{\mathcal S}\arrow{n,l}{H}\arrow{e,t}{\overline F}
\node{\mathcal S'}
\end{diagram}\]
commutes, where $H, G'$ are canonical $Ann$-equivalences, and therefore $\overline F$ induces an Ann-functor.
}
\begin{proof}
 This Proposition is naturally extended of Proposition 2 [8].
\end{proof}

{\bf Definition 3.} {\it
Let $\mathcal S=(R, M,h),\ \mathcal S'=(R',M',h')$ be   $Ann-$categories. A functor  $F: \mathcal S\rightarrow \mathcal S'$ is called a functor of the type  $(p, q)$ if 
$$F(x)=p(x),\ \ \ \ F(x,a)=(p(x), q(a)), $$ 
where \  $p:R\rightarrow R'$ is a ring homomorphism  and  $q:M\rightarrow M'$ is a group homomorphism satisfying  
$$q(xa)=p(x)q(a), q(ax)=q(a)p(x),$$ 
for $x\in R,  a\in M.$
}

{\bf Proposition 4.3.} {\it
Let $\A=(R,M,h),\  \A'=(R',M',h')$ be Ann-categories  and   $(F,\Fa,\Fm)$ is an Ann-functor  from  $\A$ to $\A'$. Then,  $(F,\Fa,\Fm)$  is a functor of the type $ ( p, q).$
}

\begin{proof}
\indent  
For $x,y\in R$, we have
  $$\Fa_{x,y}:F(x)\ts F(y)\rightarrow F(x\ts y),\Fm_{x,y}:F(x)\tx F(y)\rightarrow F(x\tx y)$$
are morphisms in the  Ann-category  $\mathcal \A'$. Hence, $F(x)+F(y)=F(x+y)$ and  
$F(x). F(y)=F(xy)$, so the map $p:R\rightarrow R'$ given  by $p(x)=F(x)$ is a ring homomorphism.

Assume  that $F(x,a)=(p(x),q_x(a))$. Since  $(F,\Fa)$ is a Gr-functor,  according to  Theorem 5 [8], $q_x=q$ for all $x\in R$. Moreover $q$ is a group homomorphism:  

\begin{equation}
q(a+b)=q(a)+q(b)\tag{2}
\end{equation}
for all $a,b\in M$.

Since  $(F,\Fm)$ is a $\tx$-functor,  the following diagram 
$$\begin{CD}
Fx\tx Fy @>\Fm >> F(x\tx y)\\
@V F((x,a))\tx F((y,b)) VV   @ VV F((x,a)\tx (y,b)) V \\
Fx\tx Fy @ >\Fm >> F(x\tx y)
\end{CD}$$
commutes, for all morphisms  $(x,a), (y,b)$. So, we have:
\begin{gather*}
F((x,a)\tx (y,b))=F(x,a)\tx F(y,b)\\
\Leftrightarrow   q_{xy}(ay+xb)=q_x(a)F(y)+F(x)q_y(b) \tag{3}
\end{gather*}
Applying $q_x=q_y=q_{xy}=q$ to the relation (3), we have:
\begin{equation}
 q(ay+xb)=q(a)F(y)+F(x)q(b) \tag{4}
 \end{equation}
Applying  $x =1$ to (4), we have:
\begin{equation}
 q(ay)=q(a)F(y)=q(a)p(y) \tag{5}
 \end{equation}
Applying  $y =1$ to (4), we have:
\begin{equation}
 q(xb)=F(x)q(b)=p(x)q(b)\tag{6}
 \end{equation}

If  $R'$-bimodule  $M'$ is regarded as an   $R$-bimodule thanks  to the actions\linebreak $xa'=p(x)a',a'x=a'p(x)$, from the equations  (2), (5), (6) we may show that  $q:M\rightarrow M'$ is a homomorphism  between  $R$-bimodules.
\end{proof}


\subsection {Classification of Ann-functors}

The existence problem of Ann-functors between  Ann-categories has been solved for the regular Ann-categories  (Theorem  5.1 [13], Theorem 4.2 [1]) thanks  to Shukla cohomology.
In this section, we will solve that problem in the  general case. 

{\bf Definition 4.} {\it 
If  $F:(R, M, h)\rightarrow (R', M', h')$ is a functor of the type  $(p, q)$, then  $F$ induces 3-cocycles  $h_{\ast}= q_{\ast}h=q(h), 
 \ {h'}^*=p^{\ast}h'=h'p$, for example
\begin{eqnarray*}
\sigma'^*(x,y,z,t)&=&\sigma'(p(x),p(y),p(z),p(t))\\
\sigma_*(x,y,z,t)&=&q(\sigma(x,y,z,t)).
\end{eqnarray*}
The function $k=p^{\ast}h'-q_{\ast}h$ is called an obstruction of the functor of the type $(p,q)$.
}
Then we have

{\bf Theorm 4.4.} {\it
The functor  $F:(R, M, h)\rightarrow (R',M', h') $ of the type $(p, q)$ is an   $Ann-$functor iff the obstruction  $\overline{k}=0$ in  $H_{MaL}^3(R, M')$.
}

\begin {proof}
Let $(F,\Fa,\Fm): (R, M,h)\rightarrow (R',M',h')$ be an  $Ann-$functor of the type $(p, q)$. Since $\Fa_{x,y}=(\bullet, \mu(x,y)),$  $\Fm_{x,y}=(\bullet, \nu(x,y))$ where $\mu, \nu:R^{2}\rightarrow M'$, we may identify $\Fa, \Fm$ with $\mu, \nu$ and call $\mu, \nu$ the {\it pair of associated functions} with $\Fa,\Fm$. According to Lemma  3.1, $(F, \Fa, \Fm)$ is compatible with the pair of constraints $(v, v')$, i.e. the diagram (1) commutes, so we have:

\begin{enumerate}
\renewcommand{\labelenumi}{$\arabic{enumi}.$}
\setcounter{enumi}{6}
\item \label{7} $\sigma'^*(x,y,z,t)-\sigma_*(x,y,z,t)=\mu(x,y)+\mu(z,t)-\mu(x+z,y+t)-\mu(x,z)\\
\indent \hspace{4.6cm}-\mu(y,t)+\mu(x+y,z+t)$
\end{enumerate}
\noindent Since  $F$ is compatible with the associativity constraint of multiplication, the  distributivity constraints of  Ann-categories $\A$ and $\A'$, we have:

\begin{enumerate}
\renewcommand{\labelenumi}{$\arabic{enumi}.$}
\setcounter{enumi}{7}
\item 
$\alpha'^*(x,y,z)-\alpha_*(x,y,z)=x\nu(y,z)-\nu(xy,z)+\nu(x,yz)-\nu(x,y)z$
\item 
$\lambda'^*(x,y,z)-\lambda_*(x,y,z)=\nu(x,y+z)-\nu(x,y)-\nu(x,z)+x\mu(y,z)-\mu(xy,xz)$
\item 
$\rho'^*(x,y,z)-\rho_*(x,y,z)=\nu(x+y,z)-\nu(x,z)-\nu(y,z)+\mu(x,y)z-\mu(xz,yz)$
\end{enumerate}
From the equations  (7)-(10), we have:
\begin{equation}
h'^*-h_*=\delta g \tag{11}
\end{equation}
where $g=(-\mu, \nu)$. Hence the obstruction of the functor $F$ vanishes in the cohomology group  $H_{MaL}^3(R, M)$.\\

Conversely, assume that the obstruction of the functor $F$  vanishes in the cohomology group  $H_{MaL}^3(R, M')$. Then there exists a 2-cochain   $g=(\mu, \nu)$ such that  $h'^*-h_*=\delta g$. Take  $\Fa, \Fm$ be functor morphisms associated with the functions $-\mu, \nu$, we can verify that     $(F,\Fa,\Fm)$ is an Ann-functor.
\end{proof}

{\bf Theorm 4.5.} {\it
If there exists an Ann-functor $(F,\Fa,\Fm): \A\rightarrow \A'$, of the type $(p, q)$ then:
\begin {enumerate}
\item  
There exists   a bijection between the set of congruence classes of Ann-functors of the type $(p,q)$ and the cohomology group  $H^2_{MaL}(R, M')$ of the ring   $R$ with coefficients in $R$-bimodule $M'$.

\item  There exists a bijection
\[Aut(F)\rightarrow Z^1_{MaL}(R, M')\]
between the group of automorphisms of the Ann-functor $F$ and the group $Z^1_{MaL}(R, M').$
\end {enumerate}
}

\begin {proof}
(a) Suppose that there exists $(F, \Fa, \Fm): \A\rightarrow \A',$ which is an Ann-functor of the type $(p,q)$. According to Theorem  4.4, we have
\[\overline{h'^*-h_*}=0.\]
Hence, there exists a 2-cochain $k$ such that 
\[h'^*-h_*=\delta k.\]
Fix 2-cochain $k$. Now, we assume that 
\[(G,\Ga,\Gm):(R, M, h)\rightarrow (R', M', h')\]
is an   Ann-functor of the type $(p,q)$. Then, from the proof of the Theorem 4.4, we have
\[h'^*-h_*=\delta g.\]
Hence, $k-g$ is a  2-cocycle. Consider the  correspondence:
\[\Phi: class(G)\mapsto class(k-g)\]
from the  set of the congruence classes of  Ann-functors of the type $(p,q)$  to the group $H^2_{MaL}(R,M')$.

First, we prove that the above correspondence is a  map.  Indeed, suppose that 
\[(G',\Ga',\Gm'):(R, M, h)\rightarrow (R', M', h')\]
is also an  Ann-functor of the type $(p,q)$ and  $u: G\rightarrow G'$ is an Ann-morphism. Since $u$ is an $\ts$-morphism as well as an $\tx$-morphism, we have:

\begin{equation}
g'=g-\delta(u) \tag{12}
\end{equation}
 So
\[k-g'=k-g+\delta(u).\]
Thus $\overline {k-g}=\overline {k-g'}\in H^2_{MaL}(R,M')$.

Now, we prove that $\Phi$ is an injection. Assume that
\[(G,\Ga,\Gm),(G',\Ga',\Gm'):(R, M, h)\rightarrow (R', M', h')\]
are  Ann-functors of the type  $(p,q)$ and satisfying 
\[\overline {k-g}=\overline {k-g'}\in H^2_{MaL}(R,M').\]
Then, there exists an 1-cochain  $u$ such that 
\[k-g=k-g'-\delta(u)\]
That means
\[g'=g-\delta(u).\]
Hence, the following diagrams:
\[\begin{diagram}
\node{G(x)\oplus G(y)}\arrow{e,t}{\Ga}\arrow{s,l}{u_x\oplus u_y}
\node{G(x\oplus y)}\arrow{s,r}{u_{x\oplus y}}
\node{G(x)\otimes G(y)}\arrow{e,t}{\Gm}\arrow{s,l}{u_x\otimes u_y}
\node{G(x\otimes y)}\arrow{s,r}{u_{x\otimes y}}\\
\node{G'(x)\oplus G'(y)}\arrow{e,t}{\Ga'}
\node{G'(x\oplus y)}
\node{G'(x)\otimes G'(y)}\arrow{e,t}{\Gm'}
\node{G'(x\otimes y)}
\end{diagram}\]
\noindent commute, it means that $u: G\rightarrow G'$ is an  $Ann$-morphism. Therefore, 
\[class (G)=class(G').\]

Finally, we must prove that the correspondence $\Phi$ is a surjection. Indeed, assume that $g$ is an arbitrary 2-cocycle. We have:
\[\delta (k-g)=\delta k-\delta g=\delta k = h'^*-h_*.\]
Then, according to Theorem 4.4, there exists an Ann-functor 
\[(G,\Ga,\Gm):(R, M, h)\rightarrow (R', M', h')\]
of the type  $(p, q)$, and the pair of isomorphisms  $\Ga, \Gm$ associated with the 2-cochain $k-g$.

Clearly, $\Phi(G)=\overline g.$
So $\Phi$ is a surjection.\\

(b) Assume  that $F=(F,\Fa,\Fm):(R, M, h)\rightarrow (R', M', h')$ is an Ann-functor of the type  $(p,q)$ and $u \in Aut(F)$. Then, from the  equation  (12)  with $g'=g$, we have  $\delta(u)=0$, i.e.,  $u\in Z^1_{MaL}(R, M')$.
\end{proof}

\section{Ann-functors and Hochschild cohomology }
 In this section, we will consider   special Ann-functors which are related  to the  low-dimensioned  Hochschild groups.

Following, we will find a condition for  the existence of  Ann-functors of the form
\[F=(F, id, \Fm):(R, M, h)\rightarrow (R', M', h')\]
of the type $(p,0)$, where $p: R\rightarrow R'$ is a ring homomorphism.

Suppose that there exists an Ann-functor
 \[F=(F, id, \Fm=\nu): (R, M, h)\rightarrow (R', M', h')\]
of the type $(p,0)$.
Then, the equations (7) \ -\  (10) turn into:
\begin{enumerate}
\renewcommand{\labelenumi}{$\arabic{enumi}.$}
\setcounter{enumi}{12}
\item 
$\sigma'^*(x,y,z,t)=0$
\item 
$\alpha'^*(x,y,z)=x\nu(y,z)-\nu(xy,z)+\nu(x, yz)-\nu(x,y)z$
\item 
$\lambda'^*(x,y,z)=\nu(x,y+z)-\nu(x,y)-\nu(x,z)$
\item 
$\rho'^*(x,y,z)=\nu(x+y,z)-\nu(x,z)-\nu(y,z)$
\end{enumerate}
and the Theorem 4.4 turns into

{\bf Corollry 5.1.} {\it
Let $p:R\rightarrow R'$ be a ring homomorphism. There exists an Ann-functor $(F, id, \Fm)$ from $(R, M, h)$ to $(R', M', h')$ of the type  $(p,0)$ iff $\overline{h'^*}=0\in H^3_{MaL}(R, M')$.
}

Each cocycle of  $\mathbb Z$-algebras due to  Hochschild is a multi-linear function. This suggests us the following definition:

{\bf Definition 5.} {\it
An Ann-functor
\[(F, id, \Fm): (R, M, h)\rightarrow (R', M', h')\]
of the type $(p,0)$ is called a {\bf strong}  \it Ann-functor if the function $\nu: R^2\rightarrow M'$ corresponding to $\Fm$ is bi-additive.
}

If $\nu$ is a normal bi-additive function then $\nu$ is a 2-cocycle of the $\mathbb Z$-algebra $R$ with coefficients in  $R$-bimodule $M'$ due to Hochschild. Then, in the equations (13)-(16), $\alpha'^*$ is a normal multi-linear function and other functions are equal to 0. So, we can identify 
\begin{equation*}
h'^*\equiv \alpha'^*=\delta (\nu),
\end{equation*}
in which $\delta(\nu)$ is an 3-coboundary of the ring $R$ with coefficients in $R$-bimodule $M'$ due to  Hochschild. Then, we have the following proposition, as a direct corrolary of 
 Theorem 4.4.

{\bf Proposition 5.2.} {\it
Let $F:(R, M, h)\rightarrow (R', M', h')$ be a functor of the type $(p,0)$. There exists a strong Ann-functor  $(F, id, \Fm)$ iff its cohomology class $\overline{h'^*}=0$ in the  cohomology group $H_{Hochs}^3(R, M')$.
}

{\bf Theorm 5.3.} {\it
If there exists an strong Ann-functor $(F, id, \Fm): (R,M,h)\rightarrow (R',M',h'),$ of the type $(p, 0)$, then:
\begin {enumerate}
\item  There exists a bijection between the set of congruence classes of strong Ann-functors of the type $(p,0)$ and the cohomology group $H_{Hochs}^2(R, M')$ of  the ring  $R$ with coefficients in  $R$-bimodule $M'$.
\item  There exists a bijection
$$Aut(F)\rightarrow Z^1_{Hochs}(R, M')$$
\noindent between the group of automorphisms of the Ann-functor $F$ and the group $Z^1_{Hochs}(R, M').$
\end {enumerate}
}

\begin {proof}
(a) The restriction $\Phi_H$ of the map $\Phi$, refered in Theorem 4.5, on the set of congruence classes of strong Ann-functors, gives us an injection to the group  $H^2_{Hochs}(R,M')$. Moreover, it is easy to see that  $\Phi_H$ is also a surjection.

\newpage

(b) Assume that  $F=(F,id,\Fm):(R, M, h)\rightarrow (R', M', h')$ is an strong $Ann-$functor of the type $(p,0)$ and  $u \in Aut(F)$. Then $u$ is bi-linear respect to the $\ts$. So  $u\in Z_{Hochs}^1(R, M')$. The converse also  holds.
\end{proof}

\begin{center}

\end{center}
\noindent{\small {\it Dept. of Mathematics, Hanoi National University of Education, Viet Nam\\
Email: nguyenquang272002@gmail.com\\
\indent \quad\ \   ddhanhdhsphn@gmail.com}}

\begin{thebibliography}{99}
\bibitem{1}  T. P. Dung, \emph{Doctoral dissertation,} Hanoi, Vietnamese, 1992.


\bibitem{2}
 
D. B. A. Epstein,   \emph{Functors between tenxored categories}, Invent. math.1, 221 -228 (1966).


\bibitem{3}

G. Hochschild, \emph{On the cohomology groups of an associative algebra,} Ann. of Math. (2) 46,  (1945). 58-67. 

\bibitem {4}

M. Jibladze and T. Pirashvili, \emph{Third Mac Lane cohomology via categorical rings,} arxiv. math. KT/0608519 v1, 21 Aug 2006.

\bibitem{5}

S. Mac Lane, \emph{Extensions and Obstruction for rings,} Illinois J. Mathematics, 2 (1958), 316-345.


\bibitem {6}

S. Mac Lane, \emph{Homologie des anneaux et des modules,} Colloque de Topologie algebrique, Louvain (1956), 55-88.
 
\bibitem {7}

N. T. Quang,  \emph{Doctoral dissertation,} Hanoi, Vietnamese, 1988.


\bibitem {8}

N. T. Quang, \emph{On Gr-functors between Gr-categories: Obstruction theory for Gr-functors of the type $(\varphi, f)$,}  arXiv: 0708.1348 v2 [math.CT] 18 Apr  2009


\bibitem {9}

N. T. Quang and P. L. Hong Anh,  \emph{On monoidal equavalences and Ann-equavalences,} arXiv: 0705.0736 v1 [math. CT] 5 May 2007.


\bibitem {10}

N. T. Quang, D. D. Hanh and N. T. Thuy, \emph{On the axiomatics of Ann-categories,} JP Journal of Algebra, Number Theory and applications, Vol 11, $No$ 1, 2008, 59 - 72.

\bibitem {11}

N. T. Quang,  \emph{Structure of Ann-categories,} arXiv. 0805. 1505 v3 [math. CT] 6 Apr 2009.

\bibitem {12}

N. T. Quang, \emph{Ann-categories and the Mac Lane-Shukla cohomology of rings,}  Abelian groups and modules $No$ 11,12 (Russian), 166 - 183, Tomsk. Gos. Univ., Tomsk, 1994  .

\bibitem {13}

N. T. Quang, \emph{Structure of Ann-categories and Mac Lane-Shukla cohomology,} East-West, J. of Math. 5 (2003), 51-66.



\bibitem{14}

U. Shukla, \emph{Cohomologic des algÐbres associativÐ,}
 Ann. Sci. Ðcole Norm. Sup. (3) 78 (1961), 163 - 209.

\end{thebibliography}
\end{document}